\documentclass[10pt, a4paper]{article}
\usepackage{amssymb}
\usepackage{mathrsfs}

\usepackage{amsmath,amssymb}
\usepackage{color}

\numberwithin{equation}{section} 

 \newtheorem{thm}{Theorem}[section]
 
 \newtheorem{lem}[thm]{Lemma}
 \newtheorem{prop}[thm]{Proposition}

 \newtheorem{conj}[thm]{Conjecture}
\newenvironment{prf}{\noindent {\it Proof} \ }{\hfill $\Box$}

\newcommand{\eqa}{\begin{eqnarray}}
\newcommand{\eeqa}{\end{eqnarray}}
\newcommand{\beq}{\begin{equation}}
\newcommand{\eeq}{\end{equation}}

\newcommand\p{\partial}

\newcommand{\nn}{\nonumber}
\newcommand{\ve}{\varepsilon}

\newcommand{\lm}{\lambda}
\newcommand{\al}{\alpha}
\newcommand{\gm}{\gamma}

\newcommand{\ta}{\theta}

\newcommand{\htv}{{\hat{v}}}
\newcommand{\hte}{{\hat{E}}}
\newcommand{\htc}{{\hat{c}}}
\newcommand{\htf}{{\hat{F}}}

\begin{document}

\title{On Symmetries of the WDVV Equations}

\author{Dingdian Xu, Youjin Zhang }
\date{Department of Mathematical Sciences, Tsinghua University\\
Beijing 100084, P. R. China} \maketitle

\begin{abstract}
For two solutions of the WDVV equations that are related by two types of symmetries of the equations given by Dubrovin,  we show that the associated
principal hierarchies of integrable systems are related by certain
reciprocal transformation, and the tau functions of the hierarchies
are either identical or related by a
Legendre transformation. We also consider relationships between the
Virasoro constraints and topological deformations of the principal hierarchies.
\end{abstract}

\section{Introduction}
The Witten-Dijkgraaf-Verlinde-Verlinde (WDVV) equations, which arises
in the study of 2d topological field theory in the
beginning of 90's of the last century,  are given by the following system of PDEs
for a function $F=F(v^1,\dots, v^n)$ :
\begin{enumerate}
\item The variable $v^1$ is specified so that
\beq\label{wdvv-1} \eta_{\al\beta}:=\frac{\p^3 F}{\p v^1\p v^\al\p
v^\beta}=\rm{constant},\quad \det(\eta_{\al\beta})\ne 0. \eeq
\item The functions $c^\al_{\beta\gm}:=\eta^{\al\nu} c_{\nu\beta\gamma}$ with
\beq\label{wdvv-2}
c_{\al\beta\gamma}=\frac{\p^3 F}{\p
v^\nu\p v^\beta\p v^\gm}, \quad
(\eta^{\al\beta})=(\eta_{\al\beta})^{-1} \eeq yield the structure
constants of an associative algebra for given $v=(v^1,\dots,v^n)$,
i.e, they satisfy \beq c_{\al\beta}^\lm
c_{\lm\gm}^\nu=c_{\gm\beta}^\lm c_{\lm\al}^\nu\quad \rm{for\ any\
fixed} \ 1\le \al, \beta, \gm, \nu\le n. \eeq
Here and in what follows summation with respect to repeated upper and lower 
indices is assumed.
\item The function $F$ satisfies the quasi-homogeneity condition
\beq\label{quasi} \p_E F=(3-d) F+\frac12 A_{\al\beta}\, v^\al
v^\beta+ B_\al v^\al+C, \eeq here the Euler vector field has the
form \beq\label{euler} E=\sum_{\al=1}^n \left((1-\frac{d}2-\mu_\al)
v^\al+r_\al\right) \frac{\p}{\p v^\al}, \eeq and $d,
A_{\al\beta}=A_{\beta\al}, B_\al, C$, $\mu_\al, r_\al$ are some
constants with $\mu_1=-\frac{d}2$.
\end{enumerate}
These equations are satisfied by the primary free energy $F$ of the matter 
sector of a 2d topological field theory with $n$ primary fields as a function of the 
coupling constants \cite{Dij2, Dij3, Witten1}. In \cite{D3, D1} Dubrovin reformulated the WDVV equations in a coordinated free form by introducing the
notion of Frobenius manifold structure on the space of the
parameters $v^1,\dots, v^n$, and revealed rich
geometric structures of the WDVV equations which are 
important in the study of several different areas of mathematical research,
such as  the theory of Gromov -
Witten invariants, singularity theory and nonlinear integrable systems \cite{D1, weyl, DZ, DZ2}. In
particular, such geometrical structures enable one to associate a
solution of the WDVV equations with a hierarchy of bihamiltonian
integrable PDEs of hydrodynamic type which is called the principal
hierarchy \cite{DZ2}. This hierarchy of integrable systems plays
important role in the procedure of reconstructing a 2D topological
field theory(TFT) from its primary free energy as a solution of the
WDVV equations. In this construction, the tau function that
corresponds to a particular solution of the principal hierarchy
serves as the genus zero partition function, and the full genera
partition function of the 2D TFT is a particular tau function of an
integrable hierarchy of evolutionary PDEs of KdV type which is
certain deformation of the principal hierarchy, such a deformation of the 
principal hierarchy is 
call the {\em{topological deformation}} \cite{DZ2}.

In this paper we are to interpret certain symmetries of the WDVV
equations in term of the associated principal hierarchy and its
tau functions. The symmetries we consider here are given by Dubrovin
in  Appendix B of \cite{D1}, where they are called symmetries of
type-1 and type-2 respectively\footnote{In \cite{CKS}
the Lie algebra of general infinitesimal symmetries of the WDVV equations without the homogeneity condition is considered.}. These symmetries are obtained from the
Schlesinger transformations of the system of linear ODEs with rational coefficients
which are associated to the Frobenius manifolds (see
Remark 4.2 of \cite{D2} for details). It turns out that in terms of the
principal hierarchies and their tau functions these symmetries have
a simple and natural interpretation. On the principal
hierarchy these symmetries act as reciprocal transformations, and
on the associated tau functions these two types of symemtries 
either keep the tau functions unchanged or act as Legendre transformations.

Recall that  a symmetry of the WDVV equations consists of
transformations 
\beq v^\al\mapsto \hat{v}^\al,\quad
\eta_{\al\beta}\mapsto\hat{\eta}_{\al\beta},\quad F\mapsto \hat{F}
\eeq 
that preserve the WDVV equations. The two types of symmetries
given in \cite{D1} have the following form:
\begin{enumerate}
\item Type-1 symmetries: they are given by the transformations defined by
\beq\label{zh-8}
{\hat v}^\al=\eta^{\al\gm}\frac{\p^2 F(v)}{\p
v^\gm\p v^\kappa},\quad {\hat\eta}_{\al\beta}=\eta_{\al\beta}, \quad
\frac{\p^2 \hat{F}(\hat{v})}{\p \htv^\al\p \htv^\beta}=  \frac{\p^2
{F(v)}}{\p v^\al\p v^\beta} 
\eeq 
for any given $1\le\kappa\le n$
such that the matrix $(c^\al_{\beta\kappa})$ is invertible. Note
that in this case the transformed function $\hat{F}(\htv)$ satisfies
the WDVV equations with $\htv^\kappa$ as the specified variable, and
the equations in \eqref{wdvv-1} is replaced by 
\beq\label{zh-1-28}
\hat{\eta}_{\al\beta}:=\frac{\p^3 \hat{F}}{\p \htv^\kappa\p
\htv^\al\p \htv^\beta}=\rm{constant},\quad
\det(\hat{\eta}_{\al\beta})\ne 0. \eeq

\item Type-2 symmetries: they are given by the transformation defined by
\begin{align}
\begin{split}\label{type-2t}
&{\hat v}^1=-\frac12\,\frac{\eta_{\sigma\gamma}  v^\sigma
v^\gamma}{v^n},\quad
{\hat v}^\al=\frac{v^\al}{v^n}\ \ \mbox{for}\ \al\ne 1, n,\quad {\hat v}^n=\frac1{v^n},\\
& {\hat\eta}_{\al\beta}=\eta_{\al\beta},\quad \hat{F}(\hat v) =
(v^n)^{-2}\left(-F(v)+\frac{1}{2} \eta_{\sigma\gamma}v^1  v^\sigma
v^\gamma\right).
\end{split}
\end{align}
 Here in this case we assume that we can normalize the coordinates $v^1,\dots, v^n$ such that the functions $\eta_{\al\beta}$ take the values
 \beq
 \eta_{\al\beta}=\delta_{\al+\beta,n+1},
 \eeq
 and in the expression of the Euler vector field $E$ given in \eqref{euler} the constants $r_\al=0$ whenever $\mu_\al\ne 1-\frac{d}2$.
We also impose the additional conditions that in the cases
 when $d=1$ and $d=2$ the constants $r_n$ and $B_1$ that appear in \eqref{euler} and \eqref{quasi} vanishes respectively. Note that the transformation \eqref{type-2t}
is obtained from the one given in Append B of \cite{D1} by changing the signs of $\htv^1, \htv^n$ and of
 $\hat{F}$, we make this modification so that the above transformation is an involution.
\end{enumerate}

We arrange the content of the paper as follows. We first recall in
Sec. \ref{sec-2} the definition of the principal hierarchy and its
tau function associated to a solution of the WDVV equations, or
equivalently, to a Frobenius manifold. In Sec. \ref{sec-3} and Sec.
\ref{sec-4} we  show respectively that the actions of the type-1 and
type-2 symmetries on the principal hierarchies are given by certain
reciprocal transformations, and we give the transformation rule
of the associated tau functions.  In Sec. \ref{sec-5} we consider the  
transformation rule of the Virasoro constraints for the tau functions of the principal hierarchy. We conclude the paper with a discussion on actions of the symmetries
of the WDVV equation on the topological deformations of the principal hierarchy.

\section{The principal hierarchy}\label{sec-2}
Given a solution of the WDVV equations, the associated principal
hierarchy consists of Hamiltonian systems of the following form:
\beq\label{zh-9} 
\frac{\p v^\al}{\p t^{\beta,q}}=\eta^{\al\gm}
\p_x\left(\frac{\p\theta_{\beta,q+1}}{\p v^\gm}\right),\quad
\al,\beta=1,\dots,n,\ q\ge 0. 
\eeq 
Here the densities $\theta_{\beta,q+1}$ of the Hamiltonians
$H_{\beta,q}=\int\theta_{\beta,q+1}(v(x)) dx$ are given by the flat
coordinates of the deformed flat connection of the corresponding
Frobenius manifold \cite{D1}. Denote 
\beq \theta_\al(z)=\sum_{p\ge 0} \theta_{\al,p}(v) z^p,\quad \al=1,\dots,n, 
\eeq 
then the functions
$\theta_{\al,p}(v)$  are determined by the equations
\begin{align}
&\label{zh-2} \p_\al\p_\beta\theta_{\nu}(z)=z\,
c^\gm_{\al\beta}\p_\gm\theta_{\nu}(z),
\quad \p_\al=\frac{\p}{\p v^\al},\quad \al, \beta, \nu=1,\dots,n.\\
&\p_E\,\p_\beta\theta_{\al,p}(v)=\left( p +\mu_\al+\mu_\beta\right)
\p_\beta\theta_{\al,p}(v)+\sum_{k=1}^p\p_\beta
\theta_{\ve,p-k}(v)\,\left(R_k\right)^\ve_\al\label{zh-4}
\end{align}
and are normalized by the conditions
\begin{align}\label{zh-7}
&\theta_{\al,0}=\eta_{\al\gm} v^\gm:=v_\al,\\
&\p_\gm\theta_\al(z) \eta^{\gm\nu}
\p_\nu\theta_{\beta}(-z)=\eta_{\al\beta}.\label{zh-3}
\end{align}
Where the  constant matrices $R_1, R_2, \dots$  have the properties
\begin{enumerate}
\item $(R_k)^\al_\beta\ne 0$ only if $\mu_\al-\mu_\beta=k$,
\item $ \eta_{\al\gamma} (R_k)^\gm_\beta=(-1)^{k+1} \eta_{\beta\gamma} 
(R_k)^\gm_\al $.
\end{enumerate}
From the first property we see that we have only finitely many
nonzero matrices $R_1,\dots, R_m$, the number $m$ is determined by
the particular solution of the WDVV equations. These matrices are
defined up to  conjugations \beq\label{cong} R=R_1+\dots+R_m\mapsto
\tilde{R}=G R G^{-1} \eeq given by nondegenerate constant matrices
$G$ that satisfy certain conditions, see \cite{D1, D2, DZ2} for
details. These matrices form part of the monodromy data
\[(V, [R], \mu, <\,,\,>, e_1)\]  of the Frobenius manifold at $z=0$.
Here $V$ is the n-dimensional vector space spanned by $e_1,\dots,
e_n$, $[R]$ is the equivalence class (w.r.t. the above conjugation)
represented by the operator $R$ that acts on $V$ by $R e_\al=
R^\gm_\al e_\gm$, the action of the operator $\mu$ on $V$ is given
by the diagonal matrix  $\mu={\mbox{diag}}(\mu_1,\dots, \mu_n)$,
and the bilinear form is given by $<e_\al,e_\beta>=\eta_{\al\beta}$.
Note that the matrix $\mu$ satsifies the anti-symmetry condition
\beq (\mu_\al+\mu_\beta) \eta_{\al\beta}=0. \eeq

The functions $\theta_{\al,p}(v)$ satisfy the following tau-symmetry
condition: 
\beq \frac{\p\theta_{\al,p}(v)}{\p t^{\beta,q}}=
\frac{\p\theta_{\beta,q}(v)}{\p t^{\al,p}}, \quad \al,
\beta=1,\dots,n. 
\eeq 
This property enables one to introduce the tau
function $\tau$ of the principal hierarchy \eqref{zh-9}. It is
defined for any given solution  $v^\al=v^\al(t)$ of the hierarchy
and required to satisfy the equations 
\beq \frac{\p^2\log\tau}{\p x\p t^{\al,p}}=\theta_{\al,p}(v(t)),\quad \al=1,\dots,n,\ p\ge 0.
\eeq 
Note that the flow $\frac{\p}{\p t^{1,0}}$ coincides with $\frac{\p}{\p x}$.
In order to fix $\tau$ up to a linear
function of $t^{\al,p}$, we are to use the functions
$\Omega_{\al,p;\beta,q}(v)$ defined by the following identities
\cite{D1}: \beq\label{zh-11} \sum \Omega_{\al,p;\beta,q}(v)\, z^p
w^q=\frac{\p_\gm\theta_\al(z)\,\eta^{\gm\xi}\p_\xi\theta_\beta(w)-\eta_{\al\beta}}{z+w}.
\eeq 
Then for any solution $v(t)=(v^1(t),\dots, v^n(t))$ of the principal
hierarchy \eqref{zh-9} we can fix, up to a factor of the
form $e^{\sum a_\al t^\al+b}$, the tau function by the following
relations: \beq\label{zh-12} \frac{\p^2\log\tau}{\p t^{\al,p}\p
t^{\beta,q}}=\Omega_{\al,p;\beta,q}(v(t)),\quad \al,
\beta=1,\dots,n,\ p, q\ge 0. \eeq

\section{Actions of the type-1 symmetries}\label{sec-3}
The main results of this section were obtained in \cite{DZ2}, we recollect them
give the proofs here for the  convenience of comparison of them with the results that are related to the type-2 symmetries of the WDVV equations.

Let $F=F(v^1,\dots, v^n)$ be a solution of the WDVV equations with
Euler vector field $E$ of the form \eqref{euler}. After the action
of a symmetry of type-1, we obtain a new solution
$\hat{F}=\hat{F}(\htv^1,\dots,\htv^n)$. Here we note that the unity
vector field is $\frac{\p}{\p \htv^\kappa}$ instead of $\frac{\p}{\p
\htv^1}$, and the first set of equations \eqref{wdvv-1} of the WDVV
equations is changed to \eqref{zh-1-28}.

We first note that in the new coordinates $\htv^1,\dots,\htv^n$ the
Euler vector field $E$ given by \eqref{euler} has the expression \beq E=\sum_{\al=1}^n
\left(\left(1-\frac{{\hat{d}}}2-\hat\mu_\al\right)
\htv^\al+\hat{r}_\al\right) \frac{\p}{\p\htv^\al}, 
\eeq 
where
\beq\label{zh-5} 
\hat{d}=-2\mu_\kappa,\quad
\hat{\mu}_\al={\mu}_\al,\quad \hat{r}_\al=A_{\kappa\xi}
\eta^{\xi\al}. 
\eeq 
Let us show that the Euler vector field
$\hat{E}$ of the new solution $\hat{F}$ coincides with $E$. In fact,
by using \eqref{quasi}, \eqref{zh-8} we have
\begin{align}
&\p_{E}\frac{\p^2\htf(\htv)}{\p\htv^\al\p\htv^\beta}=\p_{E}\frac{\p^2 F(v)}{\p v^\al\p v^\beta}=(1+\mu_\al+\mu_\beta) \frac{\p^2 F(v)}{\p v^\al\p v^\beta}+A_{\al\beta}\nn\\
&\ = (1+\mu_\al+\mu_\beta) \frac{\p^2 \htf(\htv)}{\p \htv^\al\p
\htv^\beta}+A_{\al\beta},
\end{align}
This yields the identity \beq \p_E \hat{F}(\htv)=(3-\hat{d})
\hat{F}(\htv)+\frac12 \hat{A}_{\al\beta} \htv^\al\htv^\beta+\hat{B}_\al
\htv^\al+\hat{C} \eeq for some constants $\hat{A}_{\al\beta},
\hat{B}_\al, \hat{C}$.

Now we consider the relations of the densities
$\hat{\theta}_{\al,p}(\htv)$ of the Hamiltonians of the principal
hierarchy associated to $\hat{F}(\htv)$ with the ones that are associated to
the original solution $F(v)$ of the WDVV equations.
\begin{lem}\
i)\ The functions $\hat{\theta}_{\al,p}(\htv)$ for $\htf(\htv)$ can be
determined by the relations \beq\label{zh-6}
\frac{\p\hat{\theta}_{\al,p}(\htv)}{\p \htv^\beta}=
\frac{\p\theta_{\al,p}(v)}{\p v^\beta},\quad \al,\beta=1,\dots,n,\ p\ge
0
\eeq 
and the normalization conditions
\begin{equation}\label{zh-10-12}
\hat{\theta}_{\al,0}(\htv)=\eta_{\al\gamma} \htv^\gamma,\ \al=1,\dots,n.
\end{equation}
ii)\ The monodromy data $(\hat{V}, [\hat{R}], \hat{\mu}, <\,,\,>,
\hat{e}_\kappa)$ at $z=0$ of the Frobenius manifold associated to
$\htf(\htv)$ coincide with that of the Frobenius manifold associated to
$F(v)$. Here $\hat{V}$ is the n-dimensional vector space  spanned by
$\hat{e}_1,\dots, \hat{e}_n$, $[\hat{R}]$ is the equivalence class
represented by the operator $\hat{R}$ that acts on $\hat{V}$ by
$\hat{R} \hat{e}_\al= R^\gm_\al \hat{e}_\gm$, the action of the
operator $\hat{\mu}$ on $\hat{V}$ is given by the diagonal matrix
$\hat{\mu}=\mu={\mbox{diag}}(\mu_1,\dots, \mu_n)$,  and the bilinear
form is given by $<\hat{e}_\al,\hat{e}_\beta>=\eta_{\al\beta}$.
\end{lem}
\begin{prf}
From the definition of the monodromy data \cite{D1, D2}, it follows that  in order to 
prove the lemma we only need to verify that the functions $\hat{\theta}_{\al,p}(\htv)$
determined by the conditions \eqref{zh-6} and
\eqref{zh-10-12} satisfy \eqref{zh-2}--\eqref{zh-3}.
By using
\eqref{zh-8}, \eqref{zh-2} and \eqref{zh-6} we have 
\eqa
&&\frac{\p^2\hat{\theta}_{\gm}(z)}{\p\htv^\al\p\htv^\beta}
=\frac{\p}{\p\htv^\al}\left(\frac{\p\theta_\gm(z)}{\p
v^\beta}\right) =\frac{\p^2\theta_\gm(z)}{\p v^\beta\p
v^\lambda}\frac{\p v^\lambda}{\p\htv^\al}
=z c^\nu_{\beta\lambda} \frac{\p\theta_\gm(z)}{\p v^\nu}\frac{\p v^\lambda}{\p\htv^\al}\nn\\
&&=z\htc^\nu_{\al\beta}\frac{\p\theta_\gm(z)}{\p v^\nu}
=z\htc^\nu_{\al\beta}\frac{\p\hat{\theta}_\gm(z)}{\p \htv^\nu}.\nn
\eeqa
So we prove that the functions $\hat\theta_{\al,p}(\htv)$ satisfy the recursion relation \eqref{zh-2}.
Similarly, from \eqref{zh-4} and \eqref{zh-6} it follows that
\eqa
&&\p_\hte\,\frac{\p \hat{\theta}_{\al,p}(\htv)}{\p \htv^\beta}=\p_E\,\frac{\p {\theta}_{\al,p}(v)}{\p v^\beta}\nn\\
&&=\left( p+{\mu}_\al +{\mu}_\beta\right)
\frac{\p{\theta}_{\al,p}(v)}{\p v^\beta}+\sum_{k=1}^p \frac{\p{\theta}_{\ve,p-k}(v)}{\p v^\beta}\,\left(R_k\right)^\ve_\al\nn\\
&&=\left( p+\hat{\mu}_\al +\hat{\mu}_\beta\right)
\frac{\p\hat{\theta}_{\al,p}(\htv)}{\p \htv^\beta}+\sum_{k=1}^p
\frac{\p\hat{\theta}_{\ve,p-k}(\htv)}{\p
\htv^\beta}\,\left(\hat R_k\right)^\ve_\al,
\eeqa 
which proves the validity of \eqref{zh-4} for the functions $\hat\theta_{\al,p}(\htv)$.
The relations
\eqref{zh-7} and \eqref{zh-3} hold true obviously. The lemma is
proved.
\end{prf}

Let us proceed to consider the relation between the principal hierarchies 
associated to the
solutions $F(v)$ and $\hat{F}(\htv)$ of the WDVV equations. In the principal hierarchy \eqref{zh-9} we have $\frac{\p v^\al}{\p
t^{1,0}}=\frac{\p v^\al}{\p x}$, so we may identify the time
variable $t^{1,0}$ with the spatial variable $x$ and forget the flow
$\frac{\p}{\p t^{1,0}}$ in the hierarchy. Similarly, for the
principal hierarchy
\beq\label{zh-10}
\frac{\p\htv^\al}{\p
\hat{t}^{\beta,q}}=\hat{\eta}^{\al\gm}
\frac{\p}{\p\hat{x}}\left(\frac{\p\hat{\theta}_{\beta,q+1}}{\p
\htv^\gm}\right),\quad \al, \beta=1,\dots, n,\ q\ge 0
\eeq
 we have $\frac{\p v^\al}{\p \hat{t}^{\kappa,0}}=\frac{\p v^\al}{\p \hat{x}}$,
 so we may also identify $\hat{t}^{\kappa,0}$ with the new spatial variable
$\hat{x}$ and forget the flow $\frac{\p}{\p \hat{t}^{\kappa,0}}$ in
the hierarchy. We will assume such an identification henceforth.

\begin{prop}\label{prp-zh-1}
The principal hierarchy \eqref{zh-10} associated to the solution $\htf(\htv)$ of the WDVV
equation is obtained from the principal hierarchy
\eqref{zh-9} by the following reciprocal transformation
\begin{align}\label{zh-9-b}
&\hat{x}=t^{\kappa,0},\quad \hat{t}^{1,0}=x,\\
& \hat{t}^{\al,p}=t^{\al,p},\quad (\al,p)\ne (1,0),\
(\kappa,0).\label{zh-13}
\end{align}
i.e. the principal hierarchy \eqref{zh-10} is
obtained from \eqref{zh-9} simply by exchange the role of the spatial
variable $x$ and the time variable $t^{\kappa,0}$.
Moreover,  any tau function $\tau(t)$ of the principal hierarchy
\eqref{zh-9} yields a tau function $\hat{\tau}(\hat{t})$ of  \eqref{zh-10} by the formula
\beq\label{zh-14} \hat{\tau}(\hat{t})=\tau(t)|_{t^{\kappa,0}\to
\hat{x}, x\to \hat{t}^{1,0},
t^{\al,p}\mapsto \hat{t}^{\al,p},\ (\al,p)\ne (1,0),\  (\kappa,0)}.
\eeq
\end{prop}

\begin{prf}
Assume that $(\beta,q)\ne (1,0), (\kappa,0)$. Let $v^1(t),\dots,
v^n(t)$ satisfy the principal hierarchy \eqref{zh-9}, and $\htv^\al,
\hat{\theta}_{\beta, q}(\htv)$ be defined as in \eqref{zh-8},
\eqref{zh-6}, \eqref{zh-10-12}. Then after the reciprocal
transformation \eqref{zh-9-b} we have
\begin{align}
&\frac{\p\htv^\al}{\p \hat{t}^{\beta,q}}=\frac{\p\htv^\al}{\p
t^{\beta,q}} =c^\al_{\kappa\lambda}\frac{\p v^\lambda}{\p
t^{\beta,q}}
=c^\al_{\kappa\lambda} c^{\lambda\nu}_\xi \frac{\p\theta_{\beta,q}}{\p v^\nu} v^\xi_x\nn\\
&=c^{\al\lambda}_{\kappa} c^\nu_{\xi\lambda}
\frac{\p\theta_{\beta,q}} {\p v^\nu} v^\xi_x=c^{\lambda}_{\kappa\xi}
c^{\al\nu}_{\lambda} \frac{\p\theta_{\beta,q}}{\p v^\nu} v^\xi_x
= c^{\al\nu}_{\lambda} \frac{\p\theta_{\beta,q}}{\p v^\nu} \frac{\p v^\lambda}{\p t^{\kappa,0}}\nn\\
&=\hat \eta^{\al\nu}\frac{\p}{\p v^\lambda}\left(\frac{\p\hat{\theta}_{\beta,q+1}}{\p
\htv^\nu}\right)
 \frac{\p v^\lambda}{\p t^{\kappa,0}}=\hat{\eta}^{\al\nu} \frac{\p}{\p\hat{x}}\left(\frac{\p\hat{\theta}_{\beta,q+1}}{\p \htv^\nu}\right).
\end{align}
In a similar way we can prove the validity of the above equation for
$(\beta,q)=(1,0), (\kappa,0)$. So the reciprocal transformation
\eqref{zh-9-b}, \eqref{zh-13} transforms the principal
hierarchy \eqref{zh-9} to the principal hierarchy \eqref{zh-10}.

From the definition \eqref{zh-11} of the functions
$\Omega_{\al,p;\beta,q}$ we see that the functions
$\hat{\Omega}_{\al,p;\beta,q}(\htv)$ satisfy 
\beq
\hat{\Omega}_{\al,p;\beta,q}(\htv)=\Omega_{\al,p;\beta,q}(v), \eeq
where $\htv^1,\dots, \htv^n$ are related to $v^1,\dots, v^n$ by the
relation \eqref{zh-8}. Then the relation \eqref{zh-14} follows from
\eqref{zh-12}. Thus we proved the Proposition.
\end{prf}

The principal hierarchy \eqref{zh-9} possesses a bihamiltonian
structure given by the following compatible Hamiltonian operators
\beq\label{zh-10-13} P_1^{\al\beta}=\eta^{\al\beta}\p_x,\quad
P_2^{\al\beta}=g^{\al\beta}(v) \p_x+\Gamma^{\al\beta}_\gamma(v)
v^\gamma_x. \eeq Here $(g^{\al\beta})$ is the intersection form of
the Frobenius manifold \cite{D1} associated to $F(v)$ and is defined by
\beq\label{metric2}
g^{\al\beta}(v)=E^\gamma(v) c^{\al\beta}_\gamma(v),\quad \al, \beta=1,\dots,n,
\eeq
$\Gamma^{\al\beta}_\gamma=-g^{\al\xi} \Gamma_{\xi\gamma}^\beta$
are the contravariant components of the Levi-Civita connection of
the metric $(g_{\al\beta})=(g^{\al\beta})^{-1}$.

In a similar way,  the flat metric $\hat{\eta}$ and the intersection
form $\hat{g}$ of the Frobenius manifold associated to the new
solution $\hat{F}(\htv)$ of the WDVV equations give in the same way a bihamiltonian
structure for the principal hierarchy \eqref{zh-10}. In the flat
coordinates $\htv^1,\dots, \htv^n$, the compatible Hamiltonian
operators \beq
\hat{P}_1^{\al\beta}=\hat{\eta}^{\al\beta}\p_{\hat{x}}, \quad
\hat{P}_2^{\al\beta}=\hat{g}^{\al\beta}(\htv)
\p_{\hat{x}}+\hat{\Gamma}^{\al\beta}_\gamma(\htv)
\htv^\gamma_{\hat{x}} \eeq
have the following relation with
the Hamiltonian operators \eqref{zh-10-13}:
\beq
\hat{\eta}^{\al\beta}=\eta^{\al\beta},\quad
\hat{g}^{\al\beta}(\htv)=g^{\al\beta}(v),\quad
\Gamma^{\al\beta}_\gamma(v)=c^{\xi}_{\kappa\gamma}(v)
\hat{\Gamma}^{\al\beta}_\xi(\htv). \eeq

In general a bihamiltonian system of hydrodynamic type may not be 
related to Frobenius manifold, in such cases we can still
perform the reciprocal transformation that exchanges the role of the
spatial and time variables. It was shown in \cite{XTZ} that such a transformation
preserves the bihamiltonian property of the system, 
and the transformation rule of the
bihamiltonian structure is similar to the one given above.
Such transformations are applied in \cite{DZ2} to certain bihamiltonian hierarchies
so that the transformed ones are associated to Frobenius manifolds.

Note that for a Hamiltonian system of hydrodynamic type the
transformation rule of the Hamiltonian structure under linear
reciprocal transformations was given by Pavlov in \cite{pavlov}. An
interesting problem is whether we still have similar transformation
rules when we apply linear
reciprocal transformations to a Hamiltonian or bihamiltonian system
which is certain deformation of a system of hydrodynamic type.  

\section{Actions of the type-2 symmetries}\label{sec-4}
As we did in the last section, we denote by $\hat{F}(\htv^1,\ldots,\htv^n)$ the solution of the WDVV equations that is obtained from a solution $F(v)$ by the action of the type-2 symmetry \eqref{type-2t}.  Note that in general the operator given by the gradient of the Euler vector field $\hat{E}$ for $\hat{F}(\htv)$ is non-diagonalizable, for the convenience of the  presentation of the results on the transformation rule of the principal hierarchies and their tau functions under the action of the type-2 symmetries, we assume that the in expression \eqref{euler} of the Euler vector field $E$ the constants $r_\al,\, \al=1,\dots, n$ vanish. Under this assumption the Euler vector field for $\hat{F}(\htv)$  has the expression (see Lemma B.1 of \cite{D1})
\beq
\hat{E}=\sum_{\al=1}(1-\frac{\hat{d}}{2}-\hat{\mu}_\al)\htv^\al\frac{\p}{\p\hat{v}^\al},
\eeq
where
\beq
\hat{d}=2-d,\ \hat{\mu}_1=\mu_n-1,\
\hat{\mu}_n=\mu_1+1,\ \hat{\mu}_\al=\mu_\al \ \al\neq 1,n.\label{zh-1-16a}
\eeq
In fact it coincides with the Euler vector field $E$ for the function $F(v)$.
The function $\hat{F}(\htv)$ satisfies the following quasi-homogeneity condition:
\beq \p_{\hat{E}} \hat{F}=(3-\hat{d})
\hat{F}+\frac12 \hat{A}_{\al\beta}\, \htv^\al\htv^\beta+\hat{B}_\al
\htv^\al+\hat{C}.
\eeq
Here $\hat{A}_{\al\beta},
\hat{B}_\al, \hat{C}$ are some constants.

It was shown in \cite{D1} that the functions
\[\hat\eta_{\al\beta}=\frac{\p^3\hat{F}(\htv)}{\p \htv^1\htv^\al\htv^\beta},\quad\hat{c}_{\al\beta\gamma}=\frac{\p^3 \hat{F}(\htv)}{\p \htv^\al \p\htv^\beta\htv^\gamma}\]
have the following relations with the functions $\eta_{\al\beta}, c_{\al\beta\gamma}$
defined in \eqref{wdvv-1}, \eqref{wdvv-2}:
\begin{align}
&\hat\eta_{\al\beta}=\eta_{\al\beta},\\
&\hat{c}_{\al\beta\gamma}(\htv)=-(v^n)^{-2}\frac{\p
v^\lambda}{\p \hat{v}^\al}\frac{\p v^\mu}{\p \hat{v}^\beta}\frac{\p
v^\nu}{\p \hat{v}^\gamma}c_{\lambda\mu\nu}(v). \label{identity1}
\end{align}
By taking $\gamma=1$ in \eqref{identity1} we obtain
\beq\label{identity2}
{\eta}_{\alpha\beta}=(v^n)^{-2}\,\eta_{\lambda\mu} \frac{\p v^\lambda}{\p\hat{v}^{\al}}\frac{\p v^\mu}{\p\hat{v}^{\beta}}.
\eeq
We also have the following identities which will be used below:
\beq\label{zh-16}
  -v^n\delta_{\al}^n\frac{\p v^{\mu}}{\p \hat{v}^\beta}-v^n\delta^n_\beta\frac{\p v^\mu}{\p \hat{v}^\al}
  =\frac{\p^2 v^\mu}{\p \hat{v}^\al\p \hat{v}^\beta}+\eta_{\al\beta}\delta_1^\mu v^n.
  \eeq

\begin{lem}
 i)\ The functions $\hat{\theta}_{\al,p}(\htv)$ for the solution $\hat{F}(\htv)$ of the WDVV equations are given by the following formulae:
\begin{align}\label{zh-18}
&\hat{\theta}_{1,0}(\htv)=\frac{1}{v^n},\quad \hat{\theta}_{1,p}(\htv)=(-1)^p\frac{\theta_{n,p-1}(v)}{v^n},\quad p>0,\nn \\
  &\hat{\theta}_{\al,p}(\htv)=(-1)^p\frac{\theta_{\al,p}(v)}{v^n},\quad 2\le\al\le n-1,\quad\ p\geq 0, \\
  &\hat{\theta}_{n,p}(\htv)=(-1)^{p+1}\frac{\theta_{1,p+1}(v)}{v^n},\quad p\geq 0.\nn
  \end{align}
ii)\  Let $(\hat{V}, [\hat{R}], \hat{\mu}, <\,,\,>,
\hat{e}_1)$ be the monodromy data at $z=0$ of the Frobenius manifold associated to
$\htf(\htv)$  with $\hat V$ being the linear space spanned by $\hat{e}_1,\dots, \hat{e}_n$. Then the operator $\hat{R}$ is given by the matrix elements
\beq\label{zh-15}
(\hat{R}_{k})^{\al}_{\beta}=(-1)^{k+\delta^\al_n+\delta^n_\beta}
 \left(R_{k+\delta(\al)-\delta(\beta)}\right)^{\al+(n-1)\delta(\al)}_{\beta+(n-1)\delta(\beta)}\, ,
 \eeq
here we denote by $\delta(\al)$ the difference of two Kronecker delta functions
\beq
\delta(\al):=\delta^1_\al-\delta^n_\al \label{dlt},
\eeq
and we assume that $R_l=0$ whenever $l\le 0$.
The operator $\hat\mu$ is given
by \eqref{zh-1-16a} and the bilinear form is defined by $<\hat{e}_\al, \hat{e}_\beta>=\eta_{\al\beta}$. 
\end{lem}

\begin{prf}
We need to verify that the functions $\hat{\theta}_{\al,p}(\htv)$ satisfy
the equations (\ref{zh-2})--(\ref{zh-3}). The validity of the normalization conditions
(\ref{zh-7}), (\ref{zh-3}) is easy to see from the definition \eqref{type-2t} of the new flat coordinates $\hat v^1,\dots, \hat v^n$, the identity \eqref{identity2} and
the relations
\begin{align}
&\theta_{1,1}(v)=\frac{\p F(v)}{\p v^1}=\frac12\,\eta_{\al\beta}  v^\al v^\beta,\\
& \frac{\p\theta_{\al,p}(v)}{\p v^1}=\theta_{\al,p-1}(v)+\delta_{\al,n}\delta_{p,0}
\quad \mbox{with}\ \theta_{\al,-1}(v)=0. \label{zh-1-19a}
\end{align}
The recursion relations \eqref{zh-2} for $\hat{\theta}_{\al,p}(\htv)$ can be verified by
using the identities \eqref{identity1}--\eqref{zh-16}
and \eqref{zh-1-19a}.

To prove the validity of the quasihomogeneity condition \eqref{zh-4}, let us first assume that $\al, \beta\ne 1, n$, then by using \eqref{zh-18} and \eqref{zh-15} we get
\begin{align}
 &   \p_\hte\,\frac{\p \hat{\theta}_{\al,p}(\htv)}{\p \htv^\beta}=E^\ve\frac{\p}{\p v^\ve}\left(\frac{\p v^1}{\p \hat{v}^\beta}\frac{\p}{\p v^1}+v^n\frac{\p}{\p v^\beta}\right)
    \left(\frac{(-1)^p\ta_{\al,p}(v)}{v^n}\right)\nn\\
    =&E^\ve\frac{\p}{\p v^\ve}\left(\frac{(-1)^{p+1} v^{n+1-\beta}}{v^n}\frac{\p\ta_{\al,p}(v)}{\p v^1}\right)
    +E^\ve\frac{\p}{\p v^\ve}\left((-1)^p\frac{\p\ta_{\al,p}(v)}{\p v^\beta}\right)\nn\\
    =&(-1)^{p+1}\hat{E}^{n+1-\beta}\frac{\p\ta_{\al,p}(v)}{\p v^1}
    +\frac{(-1)^{p+1}v^{n+1-\beta}}{v^n}E^\ve\frac{\p}{\p v^\ve}\left(\frac{\p \ta_{\al,p}(v)}{\p v^1}\right)\nn\\
    &+(-1)^p(p+\mu_\al+\mu_\beta)\frac{\p\ta_{\al,p}(v)}{\p v^\beta}+(-1)^p\sum_{k=1}^p\frac{\p\ta_{\ve,p-k}(v)}{\p v^\beta}(R_k)_\al^\ve\nn\\
=&(-1)^{p+1}(1-\frac{\hat{d}}{2}-\hat{\mu}_{n+1-\beta})\frac{v^{n+1-\beta}}{v^n}\frac{\p\ta_{\al,p}(v)}{\p v^1}\nn\\
    &+\frac{(-1)^{p+1}v^{n+1-\beta}}{v^n}
    \left((p+\mu_\al+\mu_1)\frac{\p\ta_{\al,p}(v)}{\p v^1}+\sum_{k=1}^p\frac{\p\ta_{\ve,p-k}(v)}{\p v^1}(R_k)_\al^\ve\right)\nn
\end{align}
\begin{align}    
    &+(-1)^p(p+\mu_\al+\mu_\beta)\frac{\p\ta_{\al,p}(v)}{\p v^\beta}+(-1)^p\sum_{k=1}^p\frac{\p\ta_{\ve,p-k}(v)}{\p v^\beta}(R_k)_\al^\ve\nn\\
    =&(p+\mu_\al+\mu_\beta)\frac{\p \hat{\ta}_{\al,p}(\htv)}{\p \hat{v}^\beta}+
    \sum_{k=1}^p\frac{\p}{\p\hat{v}^\beta}\left(\frac{(-1)^p\ta_{\ve,p-k}(v)}{v^n}\right)(R_k)_\al^\ve
    \nn\\
    =&(p+\hat{\mu}_\al+\hat{\mu}_\beta)\frac{\p \hat{\ta}_{\al,p}(\htv)}{\p \hat{v}^\beta}
    +\sum_{k=1}^p\sum_{\ve\neq 1,n}\frac{\p\hat{\ta}_{\ve,p-k}(\htv)}{\p \hat{v}^\beta}(\hat{R}_k)_\al^\ve\nn\\
  &  +\sum_{k=1}^p\frac{\p}{\p\hat{v}^\beta}\left(\frac{(-1)^p\ta_{1,p-k}(v)}{v^n}\right)(R_k)_\al^1
    +\sum_{k=1}^p\frac{\p}{\p\hat{v}^\beta}\left(\frac{(-1)^p\ta_{n,p-k}(v)}{v^n}\right)(R_k)_\al^n
    \nn\\
    =&(p+\hat{\mu}_\al+\hat{\mu}_\beta)\frac{\p \hat{\ta}_{\al,p}(\htv)}{\p \hat{v}^\beta}
    +\sum_{k=1}^p\sum_{\ve\neq 1,n}^p\frac{\p\hat{\ta}_{\ve,p-k}(\htv)}{\p \hat{v}^\beta}(\hat{R}_k)_\al^\ve\nn\\
  &  +\sum_{k=1}^{p-1}\frac{\p \hat{\ta}_{n,p-k-1}(\htv)}{\p\hat{v}^\beta}(\hat{R}_{k+1})_\al^n
    +\sum_{k=1}^p\frac{\p \hat{\ta}_{1,p-k+1}(\htv)}{\p\hat{v}^\beta}(\hat{R}_{k-1})_\al^1\nn\\
    =&(p+\hat{\mu}_\al+\hat{\mu}_\beta)\frac{\p \hat{\ta}_{\al,p}(\htv)}{\p \hat{v}^\beta}
    +\sum_{k=1}^p\frac{\p\hat{\ta}_{\ve,p-k}(\htv)}{\p \hat{v}^\beta}(\hat{R}_k)_\al^\ve.\nn
  \end{align}
The proof for the cases when $\al, \beta=1, n$ are similar. The lemma is proved.
\end{prf}

As we explained in the last section, we identify the time variable $t^{1,0}$ of the principal hierarchy \eqref{zh-9} with the spatial variable $x$. For the principal hierarchy that is associated to the solution $\htf(\htv)$ of the WDVV equations (see
\eqref{zh-20} below) we also identify the time variable $\hat{t}^{1,0}$ with the
spatial variable $\hat{x}$.

\begin{prop}\label{prp-zh-2}
The principal hierarchy
\beq\label{zh-20}
\frac{\p\htv^\al}{\p \hat{t}^{\beta,q}}
=\hat{\eta}^{\al\gm}\frac{\p}{\p\hat{x}}\left(\frac{\p\hat{\theta}_{\beta,q+1}(\htv)}
{\p\htv^\gm}\right),\quad \al,\beta=1,\dots,n,\ q\ge 0
\eeq
associated to the solution $\htf=\htf(\htv)$ of the WDVV
equation is related to the principal hierarchy
\eqref{zh-9} by the following reciprocal transformation:
\begin{align}
 &d \hat{x}=-v^n dx-\sum_{(\al,p)\ne (1,0)}\theta_{\al,p}(v)dt^{\al,p},\label{zh-24-1}\\
\begin{split}\label{zh-24-2}
&\hat t^{1,0}=\hat{x},\quad \hat{t}^{1,  p}=(-1)^p\, t^{n,p-1},\quad p\geq 1,\\
&\hat{t}^{n,p}=(-1)^{p+1} t^{1,p+1},\quad \hat{t}^{\al,p}=(-1)^p\, t^{\al,p},\quad \al\neq 1,n, \ p\ge 0.
 \end{split}
\end{align}
\end{prop}
\begin{prf}
From the definition of the reciprocal transformation  we have
\begin{align}
\begin{split}\label{zh-19}
 \frac{\p}{\p\hat x}&=
-\frac{1}{v^n}\frac{\p}{\p x},\\
  \frac{\p}{\p\hat{t}^{1,p}}&=(-1)^p \left(\frac{\p}{\p t^{n,p-1}}-\frac{\theta_{n,p-1}(v)}{v^n}\frac{\p}{\p x}\right),\quad p\geq 1,\\
  \frac{\p}{\p\hat{t}^{\al,p}}&=(-1)^p \left(\frac{\p}{\p t^{\al,p}}-\frac{\theta_{\al,p}(v)}{v^n}\frac{\p}{\p x}\right),\quad \al\ne 1,n,\ p\geq 0,\\
  \frac{\p}{\p\hat{t}^{n,p}}&=(-1)^{p+1}\left(\frac{\p}{\p t^{1,p+1}}-\frac{\theta_{1,p+1}(v)}{v^n}\frac{\p}{\p x}\right),\quad p\geq 0.
  \end{split}
\end{align}
Let $v^1(t),\ldots,v^n(t)$ be a solution of the principal hierarchy (\ref{zh-9}), and $\hat{v}^\al,\hat{\theta}_{\beta,q}(\htv)$
be defined as in (\ref{type-2t}), (\ref{zh-18}). Then  for $\al,\beta\neq 1,n$
we have
  \begin{align}
&\frac{\p\hat{v}^\al}{\p\hat{t}^{\beta,q}}
=(-1)^q\left(\frac{\p}{\p t^{\beta,q}}-\frac{\theta_{\beta,q}(v)}{v^n}\frac{\p}{\p x}\right) \left(\frac{v^\al}{v^n}\right) \nn\\
&=(-1)^q\left(\frac{\p v^\al}{\p t^{\beta,q}}\frac{1}{v^n}-\frac{\p v^n}{\p t^{\beta,q}}\frac{v^\al}{(v^n)^2}-
   \frac{\theta_{\beta,q}(v)}{(v^n)^2}\frac{\p v^\al}{\p x}+\frac{\theta_{\beta,q}(v)v^\al}{(v^n)^3}\frac{\p v^n}{\p x}\right)\nn\\
   &=(-1)^q\left(\eta^{\al\gamma} \frac{1}{v^n}\frac{\p}{\p x}\left(\frac{\p \theta_{\beta,q+1}(v)}{\p v^\gamma}\right)
   -\frac{v^\al}{(v^n)^2}\frac{\p\theta_{\beta,q}(v)}{\p x}-\frac{\theta_{\beta,q}(v)}{(v^n)^2}\frac{\p v^\al}{\p x}\right.\nn\\
   &\qquad \qquad\left.
   +\frac{v^\al \,\theta_{\beta,q}(v)}{(v^n)^3}\frac{\p v^n}{\p x}\right)\nn \\
   &=(-1)^q\eta^{\al\gamma}\frac{1}{v^n}\frac{\p}{\p x}\left(\frac{\p\theta_{\beta,q+1}(v)}{\p v^{\gamma}}
   -\frac{v_\gamma\theta_{\beta,q}(v)}{v^n}\right)\nn\\
   &=\eta^{\al\gamma}\left(-\frac{1}{v^n}\right)\frac{\p}{\p x}
   \left(\frac{1}{v^n}\frac{\p v^\varepsilon}{\p\hat{v}^\gamma}\frac{\p (-1)^{q+1}\theta_{\beta,q+1}(v)}{\p v^\varepsilon}\right)\nn\\
   &=\hat{\eta}^{\al\gamma}\frac{\p}{\p \hat{x}}\left(\frac{\p\hat{\theta}_{\beta,q+1}(\htv)}{\p\hat{v}^\gamma}\right).\nn
  \end{align}
  In a similar way we can prove the validity of the above equation for other cases of $\al,\beta$. The proposition is proved.
\end{prf}
\\

From the definition \eqref{zh-11} of the functions $\Omega_{\al,p;\beta,q}$ it follows that the functions $\hat{\Omega}_{\al,p;\beta,q}=\hat{\Omega}_{\al,p;\beta,q}(\htv)$ which are defined by the solution $\hat{F}(\htv)$ of the WDVV equations have the following expressions:
\begin{align}
\hat{\Omega}_{\al,p;\beta,q}=&(-1)^{p+q+1+\delta^n_\al+\delta^n_\beta}
\left(\Omega_{\alpha+(n-1)\delta(\al),p-\delta(\al);\beta+(n-1)\delta(\beta),q-\delta(\beta)}(v)\right.\nn\\
&\quad -\left.\htv^n \theta_{\alpha+(n-1)\delta(\al),p-\delta(\al)}(v) \,
\theta_{\beta+(n-1)\delta(\beta),q-\delta(\beta)}(v) \right),\label{relation}
\end{align}
where $\delta(\al), \delta(\beta)$ are defined in \eqref{dlt}.

\begin{lem}\label{zh-21}
The tau function of the principal hierarchy \eqref{zh-9} defined by \eqref{zh-12}
satisfies the following equations:
\begin{align}
&\frac{\partial}{\partial\hat{t}^{\alpha,p}}\frac{\partial}{\partial x}\log \tau=0,
\quad \forall (\alpha,p)\neq (1,0),\label{zh-20a}\\
&\frac{\partial}{\partial \hat{x}}\frac{\partial}{\partial x}\log\tau=-1.
\end{align}

\end{lem}
\begin{prf}
We first consider the case when $\al=1, p\ge 1$.  By using the relations
 (\ref{zh-19}) 
 we have
 \begin{eqnarray*}
    \frac{\partial}{\partial\hat{t}^{1,p}}\frac{\partial}{\partial x}\log\tau &=&
    (-1)^p(\frac{\partial}{\partial  t^{n,p-1}}\frac{\partial}{\partial x}\log\tau-
    \frac{\theta_{n,p-1}}{v^n}\frac{\partial^2}{\partial x^2}\log\tau)\\
    &=&(-1)^p(\theta_{n,p-1}-\theta_{n,p-1})=0.
  \end{eqnarray*}
Here we used the relation \eqref{zh-12} and the fact that $\Omega_{\al,p;1,0}
=\theta_{\al,p}$ which follows from \eqref{zh-11}.
For the cases when $\al\ne 1, p\ge 1$ the proof of the equation \eqref{zh-20a} is similar. When $(\al,p)=(1,0)$ we have
  \begin{eqnarray*}
    \frac{\partial}{\partial \hat{x}}\frac{\partial}{\partial x}\log\tau
    &=&
    -\frac{1}{v^n}\frac{\partial^2}{\partial x^2}\log\tau=-\frac{v^n}{v^n}=-1.
  \end{eqnarray*}
  The lemma is proved.
\end{prf}

It follows from the above lemma that up to the addition of a constant we have
\beq
\hat{x}=-\frac{\partial}{\partial x}\log\tau.\label{zh-20b}
\eeq
The constant can be absorbed by a translation of $\hat{x}$ in the definition of the reciprocal transformation \eqref{zh-24-1}, \eqref{zh-24-2}, so we will assume from now on the
validity of  \eqref{zh-20b}.
\begin{prop}\label{prp-zh-3}
  Any tau function $\tau(t)$ of the principal hierarchy (\ref{zh-9})  defined by
  \eqref{zh-12} yields a tau function $\hat{\tau}(\hat{t})$
  of the transformed principal hierarchy (\ref{zh-20}) by the following  Legendre transformation:
  \beq\label{legendre}
  \log\hat{\tau}= x\frac{\partial\log\tau}{\partial x}-\log \tau .
  \eeq
 together with the change of independent variables \eqref{zh-24-2} and \eqref{zh-20b}.
\end{prop}

\begin{prf}
We only need to prove the validity of the equation \eqref{zh-12} for $\hat\tau$ w.r.t.
its independent variables $\hat{t}^{\al,p}$.
For $\al, \beta\ne 1, n$, by using (\ref{zh-18}), (\ref{zh-19}) and \eqref{relation}
we have
  \begin{eqnarray*}
    \frac{\p^2\log\hat{\tau}}{\partial \hat{t}^{\alpha,p}\partial \hat{t}^{\beta,q}}
    &=& \frac{\p}{\p\hat{t}^{\al,p}}(-1)^q\left(\frac{\p}{\p t^{\beta,q}}-\frac{\theta_{\beta,q}}{v^n}\frac{\p}{\p x}\right)
    (-\log\tau+x\frac{\p\log\tau}{\p x})\\
    &=& (-1)^{p+q+1}\left(\frac{\p}{\p t^{\al,p}}-\frac{\theta_{\al,p}}{v^n}\frac{\p}{\p x}\right)\frac{\p\log\tau}{\p t^{\beta,q}}\\
    &=& (-1)^{p+q+1}(\Omega_{\al,p;\beta,q}-\frac{\theta_{\al,p}\theta_{\beta,q}}{v^n})=\hat{\Omega}_{\al,p;\beta,q}.
  \end{eqnarray*}
For other values of $\al, \beta$ we can verify the validity of the same equation,
so the proposition is proved.
\end{prf}

Note that the transformation \eqref{zh-24-2}, \eqref{zh-20b}, \eqref{legendre} is an involution, its inverse is given by \eqref{zh-24-2} and
\beq
x=-\frac{\p}{\p\hat{x}}\log\hat\tau,\quad \log\tau=\hat{x} \frac{\p\log\hat\tau}{\p\hat{x}}-\log\hat\tau.
\eeq

The flat metric $\hat{\eta}$ and the
intersection form $\hat{g}$ of the Frobenius manifold associated to
 $\hat{F}(\htv)$  give a
bihamiltonian structure for the principal hierarchy \eqref{zh-20}.
In the flat coordinates $\htv^1,\dots, \htv^n$, due to the
identities (\ref{identity1}), (\ref{identity2})  and  \eqref{metric2}, the compatible Hamiltonian operators
\beq
\hat{P}_1^{\al\beta}=\hat{\eta}^{\al\beta}\p_{\hat{x}}, \quad
\hat{P}_2^{\al\beta}=\hat{g}^{\al\beta}(\htv)
\p_{\hat{x}}+\hat{\Gamma}^{\al\beta}_{\gamma}(\htv)
\htv^\gamma_{\hat{x}}
\eeq
have the following relation with the Hamiltonian operators
given in \eqref{zh-10-13}:
\begin{align}
&\hat\eta_{\al\beta}\, d\htv^\al d\htv^\beta=(v^n)^{-2} \eta_{\al\beta}\, 
dv^\al dv^\beta,\nn\\
&\hat g_{\al\beta}(\htv)\, d\htv^\al d\htv^\beta=-(v^n)^{-2} g_{\al\beta}(v)\, 
dv^\al dv^\beta.\nn
\end{align}
Here $(g_{\al\beta})=(g^{\al\beta})^{-1}, (\hat g_{\al\beta})=(\hat g^{\al\beta})^{-1}$.
Such transformation rule for Hamiltonian structures of hydrodynamic
type is given in \cite{Fer}, and for more general type of Hamiltonian
structures  is recently given in \cite{SZ}.

\section{Virasoro constraints of the tau functions }\label{sec-5}
In this section, we consider the problem of how the actions of the symmetries of the WDVV equations change the Virasoro constraints of the tau functions.

It was shown in \cite{DZ3} that the principal hierarchy \eqref{zh-9} possesses an infinite number of Virasoro symmetries. In terms of its tau function these symmetries
can be represented in the form
\begin{align}
\frac{\p\log\tau}{\p s_m}=&
\sum a_m^{\al,p;\beta,q} \frac{\p\log\tau}{\p t^{\al,p}}\frac{\p\log\tau}{\p  t^{\beta,q}}+
\sum b_{m;\al,p}^{\beta,q} t^{\al,p} \frac{\p\log\tau}{\p t^{\beta,q}}\nn\\
&+\sum c_{m;\al,p;\beta,q} t^{\al,p} t^{\beta,q},\quad m\ge -1.\label{zh-24-1b}
\end{align}
The coefficients that appear in the above expressions are some constants, they
define a set of linear differential operators
\begin{align}
L_m=&\sum a_m^{\al,p;\beta,q} \frac{\p^2}{\p t^{\al,p}\p  t^{\beta,q}}+
\sum b_{m;\al,p}^{\beta,q} t^{\al,p} \frac{\p}{\p t^{\beta,q}}\nn\\
&+\sum c_{m;\al,p;\beta,q} t^{\al,p} t^{\beta,q}+\delta_{m,0}\, c
\end{align}
which give a representation of  half of the Virasoro algebra
\beq
[L_i, L_j]=(i-j) L_{i+j}+n\frac{i^3-i}{12} \delta_{i+j,0},\quad i, j \ge -1.
\eeq
Here $c$ is a constant.

For a generic solution of the principal hierarchy its tau function satisfies the Virasoro constraints \cite{DZ2}
\beq\label{zh-24-2b}
A_m(\bar{t};\tau)=0,\quad m\ge -1.
\eeq
Here we denote the r.h.s. of \eqref{zh-24-1b} by $A(t; \tau)$, and the shifted variable $\bar{t}$ in the above expression is defined by
\beq\label{zh-1-26}
\bar{t}^{\al,p}=t^{\al,p}-c^{\al,p}
\eeq
for some constants $c_{\al,p}$. In particular
in 2d topological field theory the partition functions are given by the tau functions which are specified by the Virasoro constraints with $c^{\al,p}=\delta^\al_1 \delta^p_1$ \cite{DZ3, LT}, and the first Virasoro constraint is the string equation
 \beq\label{string}
 \sum_{p\ge
1}t^{\al,p}\frac{\p\log\tau}{\p t^{\al,p-1}}+\frac12 \eta_{\al\beta}
t^{\al,0} t^{\beta,0}= \frac{\p \log\tau}{\p x}.
\eeq

Let $F(v)$ and $\hat{F}(\htv)$  be  solutions of the WDVV equations that are related by a type-1 or type-2 symmetry. We denote by
\begin{align}
\hat{A}(\hat{t};\hat\tau)=&\sum \hat{a}_m^{\al,p;\beta,q} \frac{\p\log\hat\tau}{\p \hat{t}^{\al,p}}\frac{\p\log\hat\tau}{\p  \hat{t}^{\beta,q}}+
\sum \hat b_{m;\al,p}^{\beta,q}\hat  t^{\al,p} \frac{\p\log\hat\tau}{\p \hat t^{\beta,q}}\nn\\
&+\sum \hat c_{m;\al,p;\beta,q} \hat t^{\al,p} \hat t^{\beta,q},\quad m\ge -1
\end{align}
the r.h.s. of the Virasoro symmetries of the principal hierarchy \eqref{zh-10} or \eqref{zh-20} associated to $\hat{F}(\htv)$.

\begin{prop}
Let $\tau(t)$ be the tau function of the principal hierarchy \eqref{zh-9} associated to the solution $F(v)$ of the WDVV equations.
Then for the tau function $\hat\tau(\hat{t})$ obtained from $\tau(t)$ by applying the type-1 symmetry we have
\beq\label{zh-1-26-2}
\hat{A}_m(\hat{t};\hat\tau)= A_m(t;\tau),\quad m\ge -1,
\eeq
and  for the tau function $\hat\tau(\hat{t})$ obtained from $\tau(t)$ by applying the type-2 symmetry we have
\beq\label{zh-1-26-3}
\hat{A}_m(\hat{t};\hat\tau)=(-1)^{m+1} A_m(t;\tau),\quad m\ge -1.
\eeq
\end{prop}
\begin{prf}
The validity of \eqref{zh-1-26-2} follows from Proposition \ref{prp-zh-1} obviously. To verify the validity of \eqref{zh-1-26-3}, we note that
the relations \eqref{zh-19} and \eqref{legendre} yield
  \begin{align}
    \frac{\p \log\hat{\tau}}{\p\hat{t}^{\al,p}}
    =&(-1)^p\left(\frac{\p}{\p t^{\al,p}}-\frac{\theta_{\al,p}(v)}{v^n}\frac{\p}{\p x}\right)
       \left(-\log\tau+x\frac{\p\log\tau}{\p x}\right)\nn\\
    =&(-1)^{p+1}\frac{\p\log\tau}{\p t^{\al,p}}+(-1)^p\left(x\frac{\p^2\log\tau}{\p t^{\al,p}\p x}
    -\frac{x\ta_{\al,p}(v)}{v^n}\frac{\p^2 log\tau}{\p x^2}\right)\nn\\
    =&(-1)^{p+1}\frac{\p\log\tau}{\p t^{\al,p}}\nn
  \end{align}
for $\al\neq 1,n, p\geq 0$, and
\begin{align}\label{zh-28}
&\frac{\p \log\hat{\tau}}{\p\hat{t}^{1,0}}=-t^{1,0},\quad  \frac{\p\log\hat{\tau}}{\p\hat{t}^{1,p}}=(-1)^{p+1}\frac{\p\log\tau}{\p t^{n,p-1}},\ p\geq 1,\nn\\
& \frac{\p\log\hat{\tau}}{\p\hat{t}^{n,p}}=(-1)^p\frac{\p\log\tau}{\p t^{1,p+1}},\ p\geq 0.\nn
\end{align}
From the above equations we obtain
\begin{align}
\hat{A}_{-1}(\hat{t};\hat\tau)=&
     \sum_{p\geq 1}\hat{t}^{\al,p}\frac{\p \log \hat{\tau}}{\p \hat{t}^{\al,p-1}}
     +\frac{1}{2}\eta_{\al\beta} \hat{t}^{\al,0}\hat{t}^{\beta,0}\nn\\
   =&\sum_{p\geq 1,\al\neq 1,n}\hat{t}^{\al,p}\frac{\p\log\hat\tau}{\p \hat{t}^{\al,p-1}}
     +\sum_{p\geq 1}\hat{t}^{1,p}\frac{\p\log\hat{\tau}}{\p \hat{t}^{1,p-1}}
     +\sum_{p\geq 1}\hat{t}^{n,p}\frac{\p\log\hat{\tau}}{\p \hat{t}^{n,p-1}}\nn\\
     &+\sum_{\al\neq 1,n}\frac{1}{2}\eta_{\al\beta} \hat{t}^{\al,0}\hat{t}^{\beta,0}
     +\hat{t}^{1,0}\hat{t}^{n,0}\nn\\
   =&\sum_{p\geq 1,\al\neq 1,n}t^{\al,p}\frac{\p\log\tau}{\p t^{\al,p-1}}+\sum_{p\geq 2}t^{n,p-1}\frac{\p\log\tau}{\p t^{n,p-2}}
   +\sum_{p\geq 1}t^{1,p+1}\frac{\p\log\tau}{\p t^{1,p}}\nn\\
   &+\hat{t}^{1,1}\frac{\p\log\hat{\tau}}{\p \hat{t}^{1,0}}+\frac{1}{2}\sum_{\al\neq 1,n}\eta_{\al\beta} t^{\al,0}t^{\beta,0}
   +t^{1,1}\frac{\p\log\tau}{\p t^{1,0}}\nn\\
   =&\sum_{p\geq 1}t^{\al,p}\frac{\p\log\tau}{\p t^{\al,p-1}}+\frac{1}{2}\eta_{\al\beta}t^{\al,0}t^{\beta,0}
   =A_{-1}(t;\tau).\nn
\end{align}
The proof of the relation \eqref{zh-1-26-3} for $m\ge 0$ is similar, so we omit it here. The proposition is proved.
\end{prf}

From the above proposition we see that after the action of the tyep-1 and type-2 symmetries of the WDVV equations,  the topological solution of the principal
hierarchy \eqref{zh-9} that is specified by the Virasoro constraints \eqref{zh-24-2b},
\eqref{zh-1-26} with $c^{\al,p}=\delta^\al_1 \delta^p_1$ is transformed to a tau function of the principal hierarchies \eqref{zh-10} and \eqref{zh-20} respectively, they satisfy the  Virasoro constraints
\beq\label{zh-1-29}
\hat A_m(\bar{\hat t};\hat\tau)=0,\quad \bar{\hat t}^{\al,p}=\hat t^{\al,p}-\hat{c}^{\al,p}.
\eeq
with
\begin{align}
\hat c^{\al,p}=\left\{\begin{array}{ll} \delta^{\al}_{1}\delta^p_1,\quad &\mbox{for the type-1 symmetry,}\\
-\delta^\al_n\delta^p_0,\quad &\mbox{for the type-2 symmetry.}
\end{array}\right.
\end{align}
Note that the tau function $\hat\tau(\hat t)$ for the topological solution of the principal hierarchy \eqref{zh-10} satisfies the Virasoro constraints  \eqref{zh-1-29}
with $\hat c^{\al,p}=\delta^\al_\kappa \delta^p_1$.

\section{Conclusion}\label{sec-6}
For two solutions of the WDVV equations related by the type-1 or type-2 symmetries, we have shown that the associated principal hierarchies are related by certain reciprocal transformation, and their tau functions are either identical or related by a Legendre transformation. We also considered the relation of the Virasoro constraints for their tau functions.

It was shown in \cite{DZ2} that the principal hierarchy associated to a semisimple Frobenius manifold has a unique deformation of the form
\beq\label{zh-1-27}
\frac{\p w^\al}{\p t^{\beta,q}}= \eta^{\al\gm}
\p_x\left(\frac{\p\theta_{\beta,q+1}(w)}{\p w^\gm}\right)+
\sum_{g\ge 1} \ve^{2 g} K^\al_{\beta,q;g}(w;w_x,\dots,w^{(2g+1)}).
\eeq
Here the $K^\al_{\beta,q;g}$ are polynomials of $w^\gamma_x,\dots,
\p_x^{2g+1} w^\gamma$ with coefficients depending smoothly on $w^1,\dots, w^n$.
Such a deformation is called the topological deformation of the principal hierarchy.
It preserves the tau structure of the principal hierarchy and has an infinite number of Virasoro symmetries. Moreover, in terms of the tau function the Virasoro symmetries are required to be linearized, i.e. they can be represented by
\[\frac{\p\tau(t;\ve)}{\p s_m}=\ve^2 L_m|_{t^{\al,p}\to \ve t^{\al,p}} \tau(t;\ve),\quad m\ge -1.\]
Here and in what follows we use $\tau(t;\ve)$ to denote the tau function of the topological deformation of the principal hierarchy, and we redenote by
\[\tau^{[0]}(x,t)=e^{{\cal F}_0(x,t)}\]
the tau function of the principal hierarchy.
For a semisimple Frobenius manifold that is associated to a 2d topological field theory, the topological deformation of the principal hierarchy is supposed to determine the partition function of the model via its tau function specified by the Virasoro constraints
\beq
L_m|_{t^{\al,p}\to \ve t^{\al,p}-\delta^\al_1\delta^p_1}\,\tau(t;\ve)=0,\quad m\ge -1.
\eeq
When $m=-1$ the above constraint is just the string equation \eqref{string}.

We then  have the following natural question: For any two semisimple Frobenius manifolds related by the type-1 or type-2 symmetries of the WDVV equations, what is the relationship between the topological deformations of their principal hierarchies?

From the construction of the topological deformation of the principal hierarchy given in \cite{DZ2} we know that
the deformed hierarchy \eqref{zh-1-27} is related to the principal hierarchy \eqref{zh-9} via a so called qusi-Miura transformation of the form
\beq
w^\al=v^\al+ \eta^{\al\gamma}
\frac{\p^2}{\p x\p t^{\gamma,0}} \sum_{g\ge 1}\ve^{2g} F_g(v;v_x,\dots,\p_x^{3g-2}v),\quad \al=1,\dots,n.
\eeq
Here the functions $F_g$ are determined by the loop equation associated to the semisimple Frobenius manifold \cite{DZ2}. In particular, we have
\beq\label{f-tp}
F_1(v,v_x)=\frac1{24}\det{\left(c_{\al\beta\gamma}(v) v^\gamma_x\right)}+G(v),
\eeq
where $G(v)$ is the G-function of the Frobenius manifold \cite{DZ}.
The tau function $\tau(t;\ve)$ of the deformed hierarchy
is related to a solution of the deformed hierarchy by the formula
\[w^\al(t)=\ve^2 \eta^{\al\gamma}\frac{\p^2\log\tau(t;\ve)}{\p x\p t^{\gamma,0}},\]
and the tau function has the genus expension
\beq\label{gep-1}
\tau(t;\ve)=e^{\sum_{g\ge 0} \ve^{2g-2} {\cal{F}}_g(t)}.
\eeq
Here ${\cal{F}}_g(t)=F_g(v(x,t),\dots, \p_x^{3g-2} v(x,t))$.

For the type-1 symmetry of the WDVV equations, we see from the above mentioned construction that the topological deformation of the principal hierarchy associated to $\hat{F}(\htv)$ is obtained from that of the principal hierarchy \eqref{zh-9} by using $t^{\kappa,0}$ as the new spatial variable $\hat{x}$.
The time variables are given by $\hat t^{\al,p}=t^{\al,p}$ for $ (\al,p)\ne (\kappa,0)$ and $\hat{t}^{\kappa,0}=\hat{x}$, and the tau functions of these topological deformations of the principal hierarchies are related by $\hat\tau(\hat x,\hat t)=\tau(x,t)$.

Based on the results of  Propositions \ref{prp-zh-2}, \ref{prp-zh-3} we have the following conjecture for the type-2 symmetry of the WDVV equations.
\begin{conj}
For the type-2 symmetry of the WDVV equations, the topological deformation of the principal hierarchy \eqref{zh-20} associated to $\hat{F}(\htv)$ is obtained, up to a Miura type transformation, from that of the principal hierarchy \eqref{zh-9} by the following Legendre transformation of the tau function:
\begin{align}\label{transform2}
 \begin{split}
  &\log\hat{\tau}(\hat{t};\ve)=\log\tau(t;i\ve)-x\frac{\p\log\tau(t;i\ve)}{\p x}\\
  &\hat{x}=\ve^2\frac{\p\log{\tau}(t;i\ve)}{\p x},\quad \hat{t}^{1,p}=(-1)^{p}t^{n,p-1},\quad p\ge 1\\
  &\hat{t}^{n,p}=(-1)^{p+1}t^{1,p+1},\quad \hat{t}^{\al,p}=(-1)^p t^{\al,p}, \quad \al\neq 1,n,\ p\ge 0.
  \end{split}
\end{align}
\end{conj}

Let us explain the validity of this conjecture at the approximation up to $\ve^2$. To this end we perform a genus expansion of $\hat\tau(\hat t;\ve)$ as follows:
\beq
\hat\tau(\hat t;\ve)=e^{\sum_{g\ge 0} \ve^{2g-2} \hat{\cal{F}}_g(\hat t)}.
\eeq
Then by using the genus expansion \eqref{gep-1} of the tau function $\tau(t;\ve)$ we can rewrite the first equation of \eqref{transform2} in the form
\begin{align}
&\ve^{-2} \left(\hat{\cal F}_0(\hat t_0)+\ve^2 \frac{\hat{\cal F}_0(\hat t_0)}{\p \hat x_0}\frac{\p{\cal F}_1(t)}{\p x} \right)+\hat{\cal F}_1(\hat t_0)+{\cal O}(\ve^2)\nn\\
&=-\ve^2 {\cal F}_0(t)+{\cal F}_1(t)-x\left(-\ve^{-2}\frac{\p{\cal F}_0(t)}{\p x}+
\frac{\p{\cal F}_1(t)}{\p x}\right)+{\cal O}(\ve^2) .\label{exp-1}
\end{align}
Here we expand $\hat t^{\al,p}=\hat t^{\al,p}(t;\ve)$ that are defined in \eqref{transform2} in the form
\[\hat t^{1,0}=\hat x=\hat x_0+\ve^2\hat x_1+{\cal O}(\ve^4)=-\frac{\p {\cal F}_0(t)}{\p x}+\ve^2 \frac{\p{\cal F}_1(t)}{\p x}+{\cal O}(\ve^4),\]
and $\hat t^{\al,p}=\hat t^{\al,p}_0$ for $(\al,p)\ne (1,0)$.
By comparing the coefficients of $\ve^{-2}$ of the left and right sides of \eqref{exp-1} we get
\beq
\hat{\mathcal{F}}_0(\hat{t}_0)=-\mathcal{F}_0(t)
  +x\frac{\p \mathcal{F}_0(t)}{\p x}.\label{freeeng0}
  \eeq
From this it follows that 
\[x=-\frac{\p \hat {\cal F}_0(\hat t_0)}{\p \hat x_0}.\]
Then the coefficients of $\ve^0$ of the equation \eqref{exp-1} yields
\beq
\hat{\mathcal{F}}_1(\hat{t}_0)=\mathcal{F}_1(t).\label{freeeng1}
\eeq

The formula \eqref{freeeng0} coincides with the Legendre transformation \eqref{legendre} between the tau functions of the
principal hierarchies. From the formula \eqref{freeeng1} we get
\beq
\hat{\cal{F}}_1(\hat t_0)=F_1(v,v_x)|_{v=v(t)}=(\frac{1}{24} \det\left(c_{\al\beta\gamma}(v) v^\gamma_x\right)+G(v))|_{v=v(t)}.
\eeq
On the other hand, from \eqref{identity1} it follows that
\beq
\frac{1}{24}\log \det(\hat{c}_{\al\beta\gm}(\hat{v})\hat{v}_{\hat{x}_0}^\gm)
  =\frac{1}{24}\log \det(c_{\al\beta\gm}(v)v_{x}^{\gm})-\frac{n}{24}\log v^n.
\eeq
And by using the result of \cite{STR} the G-function of the Frobenius manifold associated to $\hat{F}(\htv)$ is given by
\beq
\hat{G}(\hat{v})=G(v)+(\frac{n}{24}-\frac{1}{2})\log v^n.
\eeq
Thus we arrive at
\begin{align}
\hat{\cal{F}}_1(\hat t_0)&=(\frac{1}{24}\log \det(\hat{c}_{\al\beta\gm}(\hat{v})\hat{v}_{\hat{x}_0}^\gm)+\hat{G}(\htv)+\frac12 \log v^n)|_{v=v(t)}\nn\\
&=(\frac{1}{24}\log \det(\hat{c}_{\al\beta\gm}(\hat{v})\hat{v}_{\hat{x}_0}^\gm)+\hat{G}(\htv)-\frac12 \log \htv^n)|_{\htv=\htv(\hat t_0)}.
\end{align}
Then from the formula \eqref{f-tp} for the topological deformation of the principal hierarchy and the
equations
\beq
\frac{\p\hat{v}^n}{\p \hat{t}^{\gamma,0}}=\eta_{\gamma\sigma} \htv^\sigma_{\hat x},\quad \al=1,\dots,n
\eeq
we see that the hierarchy of equations satisfied by
\beq
\hat w^\al=\ve^2 \eta^{\al\gamma} \frac{\p^2 \log\hat\tau(\hat t;\ve)}{\p\hat x \p \hat t^{\gamma,0}},\quad \al=1,\dots,n
\eeq
is related, at the approximation up to $\ve^2$, to the topological deformation of the principal hierarchy associated to $\hat{F}(\htv)$ by the following Miura type transformation
\beq
\hat w^\al\mapsto \hat w^\al+
\ve^2 \frac{{\hat w}^\al_{\hat x \hat x} {\hat w}^n-{\hat w}^\al_{\hat x}{\hat w}^n_{\hat x}}{2 (\hat w^n)^2}+{\cal O}(\ve^4),\quad \al=1,\dots, n.
\eeq
We will return to the extension of the Legendre transformation to the topological deformation of a principal hierarchy in seperate publications.

\vskip 0.5truecm
\noindent{\bf Acknowledgments.}
The authors thank Boris Dubrovin and Si-Qi Liu for helpful discussions and comments.
The work is partially supported by the National Basic Research Program of China (973 Program)  No.2007CB814800 and the NSFC No.10631050.


\end{document}